\documentclass[12pt]{article}%
\usepackage{amsmath,amsfonts,amssymb,amsthm}
\usepackage{amsmath}
\usepackage{amsfonts}
\usepackage{amssymb}
\usepackage{graphicx}%
\begin{document}

\begin{center}
{\Large  \bf The mathematics of
Donald Gordon Higman}
\vskip  0.8cm

version 14 December, 2008

 \vskip 1.5cm

\begin{center}

Eiichi Bannai

Department of Mathematics

Kyushu University

Fukuoka, Japan

bannai@math.kyushu-u.ac.jp

\bigskip

Robert L.~Griess, Jr.
\\[0pt]
Department of Mathematics\\[0pt] University of Michigan\\[0pt]
Ann Arbor, MI 48109  USA \\[0pt]
 rlg@umich.edu

\bigskip

Cheryl E. Praeger

School of Mathematics and Statistics (M019)

University of Western Australia

35 Stirling Highway

Crawley WA 6009, Australia

praeger@maths.uwa.edu.au

\bigskip

Leonard Scott

Department of Mathematics

University of Virginia

Charlottesville, Virginia  USA

lls2l@virginia.edu

\end{center}

\end{center}

\newpage

\tableofcontents

\newtheorem{thm}{Theorem}%
[section]
\newtheorem{prop}[thm]{Proposition}
\newtheorem{lem}[thm]{Lemma}
\newtheorem{rem}[thm]{Remark}
\newtheorem{coro}[thm]{Corollary}
\newtheorem{conj}[thm]{Conjecture}
\newtheorem{de}[thm]{Definition}
\newtheorem{hyp}[thm]{Hypothesis}%

\newtheorem{nota}[thm]{Notation}
\newtheorem{ex}[thm]{Example}
\newtheorem{proc}[thm]{Procedure}%

\section{Introduction}

Donald Gordon Higman, (born
September 20, 1928, Vancouver, B.C., Canada),
an architect of important theories in finite groups,
representation theory, algebraic combinatorics and geometry, and a
longtime faculty
member at the University of Michigan (1960-1998),
died after a long illness on
13 February, 2006.

Don left a significant legacy of mathematical work and
personal impact on many mathematicians.
A committee was formed in 2006 to work with the Michigan Mathematical Journal
and create a memorial. The contributors have some mathematical closeness to
Don. Several of Don's fifteen doctoral students are included in this
group. The breadth of topics and quality of the writing is impressive.
For example, the
article of Brou\'e is especially direct in examining the impact of one of
Higman's basic results in representation theory (the ``Higman criterion'').

Don Higman was a serious intellectual who had the manner of a kind uncle or
concerned friend. He worked broadly in algebra and combinatorics.
He thought deeply about the ideas in his mathematical sphere, and
his style was to seek  the essence of a theory.
His work had great
influence on future developments.   This is exemplified by one of his
theorems in permutation groups, as related by Peter Neumann:
Don's ``fundamental observation that a permutation group is
primitive if and only if all its non-trivial orbital graphs are
connected changed the character of permutation group theory.
 It's a simple thing, but it introduces a point of view that allowed
lovely simplifications and extensions of the proofs of many classical
theorems due to Jordan, Manning and Wielandt.''

Len Scott relates Don's reaction to a John Thompson lecture, around 1968,
at a conference at the University of Illinois. This was not long after the
discovery of the {\it Higman--Sims sporadic simple group}. Thompson expressed agreement with
Jacques Tits's  ``heliocentric view of the universe, with the general linear group
as the sun, and these sporadic groups as just asteroids." Len happened to be
on the same elevator with Don, shortly after the lecture, when one of the
participants asked Don what he thought of the heliocentric model. Don's reply was, ``Well, it hurts your eyes to look at the sun all
the time."

The elevator passengers had a good laugh, and it really was a
marvelous line. But, reflecting further, not only can we see a part of Don's
personality and humor here, but also some of his identity as a mathematician,
and even some of his place in mathematical history.

   The decade which followed saw a dramatically changed picture of the finite simple groups. At the beginning of this
period, there were a number of well-defined infinite families of finite simple groups which contained
all known finite simple groups,
with only finitely many exceptions.  At that time, these finitely many  finite simple groups that did not fit were called
``sporadic" (group theorists were unsure as to whether more infinite families might be found).   The advance of the
classification during that period suggested that there were, indeed, no more infinite families and there were only finitely
many sporadic groups.   
Nevertheless, during the same
period, many  new sporadic groups were discovered, bringing the number of sporadics to 26 by 1975. In addition, these discoveries inspired
new algebraic constructions and conjectures, as well as striking empirical observations,
 starting in the late 1970s,
linking sporadic group theory
with other areas of mathematics and physics (this body of phenomenon is generally called ``moonshine'').   These sporadic groups may well have as much impact outside group
theory as within it.  Don Higman's independent and skeptical opinions seem to have been vindicated by the end of the 1970s.

\subsection{Acknowledgements}

The editors thank
the following individuals for their correspondence or help with
writing this biography, including contributions of photos and unpublished writing and tables:

Jonathan Alperin, Michael Aschbacher, Rosemary Bailey, Francis  Buekenhout, Peter Cameron, Charles Curtis,
David Foulser,
George Glauberman, David Gluck,
James A. (Sandy) Green, Karl Gruenberg, Willem Haemers, Jonathan Hall, Betty Higman, Sylvia Hobart, Alexander Ivanov, Misha Klin, Ching-Hung Lam, Robert Liebler, Hirosi Nagao, Peter Neumann, Manley Perkel, Geoffrey Robinson, Betty Salzberg, Alyssa Sankey, Ernie Shult, Charles Sims, Don Taylor, John Thompson, Andrew Vince, Paul-Hermann Zieschang.

Robert Griess is grateful to the University of Michigan, National Cheng Kung University, Zhejiang University, and the National Science Foundation  for financial support during the writing of this biography.

Cheryl Praeger is grateful for the support of a Federation Fellowship of the Australian Research Council.

Len Scott thanks the National Science Foundation for its support during the writing of this biography.

\section{A brief overview of Donald Higman's impact on mathematics.}

Donald Higman's beginnings were in group theory and representation theory,
clearly reflecting the interests of his Ph.\thinspace D.  advisor Reinhold Baer. Don's focal
subgroup theorem is logically at the center of basic local group theory
(the study of finite groups by normalizers of their nonidentity $p$-subgroups, where $p$ is
a prime).
The {\it focal subgroup  theorem}  has been desribed by Jon Alperin
\cite{alperin} as ``a basic concept in the theory of transfer".  Alperin also writes, ``[Don] introduced {\it relatively projective
modules}, which are the key tool for Green's whole theory of modules, and
contributed to relative homological algebra.  He proved, as a result, that if $F$
is a field of prime characteristic $p$ and $G$ is a finite group then there are
finitely many isomophism classes of indecomposable $FG$-modules if, and only if,
the Sylow $p$-subgroups of $G$ are cyclic."

Don's interest moved to groups and geometry, establishing with McLaughlin a
characterization of rank 2 linear groups. This theory was later subsumed in
the theory of BN-pairs developed by Jacques Tits. His thinking became more
abstractly combinatorial. His theory of rank 3 groups assisted in the
discovery and construction of several sporadic simple groups. Indeed, Charles
Sims and Don Higman used this theory in 1967 to make a rank 3 extension of the
Mathieu group $M_{22}$ to a then-new sporadic simple group, the Higman--Sims group.

A golden era followed, full of combinatorial theories and applications. One finds {\it coherent configurations,
association schemes}, and other theories.
Both combinatorial theorists and group theorists---especially those interested in permutation groups---realized
that algebras each had been studying had a common framework, which could be extended and studied further
with the insights provided by each discipline. The body of work created at this time continues to have an impact
and interest today.

Don Higman's work on permutation groups
transformed that  area by introducing
combinatorial techniques. Major contributions include his theory of rank 3 permutation groups and its
impact on the theory of strongly regular graphs, his introduction
of intersection matrices for permutation groups that led among other things to the study of
{\it distance transitive graphs}, and his theory of coherent configurations.
The influence this had on young researchers at the time cannot be over-stated, for example,
Peter Cameron refers to this as ``the most
exciting development'' during his time as a research student, adding that,
whereas ``Sims's methods were graph-theoretic, based
on the work of Tutte, Higman's were more representation-theoretic, and
grew from the work of Wielandt''; Michael Aschbacher regards Higman's rank 3
theory as a ``central tool'' in his early research; to Cheryl Praeger, Higman's papers on rank 3
groups and intersection matrices seemed ``like magic'', and were ``more dissected than read'' by graduate students and faculty at Oxford.
According to Eiichi Bannai,  Don Higman ``opened a
new path'' which led to the development of the area now called
algebraic
combinatorics, and is regarded as one of its founders.

The class of rank 3 permutation groups contains many important families, including
many classical groups acting on isotropic points of a bilinear form or singular points of a quadratic form.
Higman's rank 3 theory was exploited by several mathematicians to find or analyze new rank 3 groups.
In particular, Higman's table [DGH57] of numerical and group theoretic information, of actual rank 3 groups
as well as ``possibilities'', was recalled by Francis Buekenhout as
``a document of great historical value'', covering around 400 groups. Indeed both the table and the rank 3 theory
contributed to the triumph   of the collaborators Charles Sims and Donald Higman in
discovering and constructing, on 3 September 1967, the Higman Sims group. The discovery of this sporadic simple group, jointly with Charles Sims,
is one of the most famous results of Donald Higman, (see Section 5.6 
and also a contemporary account in \cite{BrSa}, and [DGH20], a preprint of which was handed out at the summer 1968 group theory conference in Ann Arbor).

\section{ Education, influences and employment}

Don Higman attended college at the University of British Columbia.  Betty Higman explains that,
``there was a woman at U.B.C. -- Celia Kreiger who urged
[Don] to continue his interest in Math and apply to grad school. She might have
been instrumental in suggesting Reinhold Baer and the University of Illinois.
Then there was Hans Zassenhaus [with whom he worked] at McGill.''
A group photo of a Canadian summer mathematics
program, probably 1951 or 1952, includes Don Higman and Hans Zassenhaus.
After Don's doctorate in 1952, he spent two years as a National Research Council
Fellow at McGill University, then two years on the faculty of Montana State
University. Thereafter, he accepted an assistant professorship at the
University of Michigan and became professor in 1960.

The main mathematical influences on Donald Higman seemed to be his advisor
Reinhold Baer and algebraist Jack E. McLaughlin, a colleague at the University
of Michigan who was older by a few years.  There was another significant contact during graduate school:
Michio Suzuki and Don Higman overlapped at the University of Illinois for
about a year.  It is clear from their later work that mutual influence was likely.
Both studied with Reinhold Baer. Coincidentally, Eiichi Bannai, an editor of
this special issue, thought of Michio Suzuki as a mathematical father and Don Higman as a
kindly mathematical uncle.

McLaughlin and Higman authored three joint articles
[DGH10], [DGH13],  [DGH16]
and had ongoing dialogue
about mathematics for decades. Their common interests included general algebra, group theory, representation theory and cohomology.

In 1968,
McLaughlin discovered and constructed a sporadic simple group using the rank 3
theory of Higman and the embedding of $PSL(3,4)$ in $PSU(4,3)$ discovered by
H. H. Mitchell \cite{mitchell}. See the book of Brauer--Sah \cite{BrSa} for how finite
simple group theory looked in the late 60s.

\section{Details of early work: group theory and representation theory}

Don
Higman's focal subgroup theorem is an insightful result about
intersections of normal subgroups of a finite group with a Sylow $p$-subgroup.
In early work on homological aspects of group representation theory, Don
established the important concept of a relatively projective module and gave a
criterion, which bears his name, for relative projectivity. Higman proved that
the finiteness of the number of isomorphism types of indecomposable modules in characteristic $p$ for a finite group
$G$ is equivalent to cyclicity of the Sylow $p$-group of $G$. This is a basic result
in the theory of modules for group algebras. He did some of the earliest
computations of degree 1 cohomology of classical groups as part of his study
of flag-transitive groups.

\subsection{Focal subgroup theorem}

Higman's theory of the focal subgroup of a Sylow subgroup was a basic tool in
local analysis in finite group theory. It could be viewed as a contribution
along the lines of Burnside, Frobenius, Gr\"un, et al. to the determination of
quotients which are $p$-groups ($p$ is a prime number). It turns out that the
focal subgroup theory is logically at the center of $p$-local group theory.
Consider the recent words of George Glauberman: ``When I was preparing some
lectures on local analysis \cite{gg}, I had to omit the
details of the definition of the transfer mapping and all proofs using the
explicit definition. I found that I could prove all the applications of
transfer that I needed by assuming without proof only two applications: the
Focal Subgroup Theorem and the fact that for a Sylow subgroup $S$ of a group
$G$, $S \cap Z(G) \cap G{}^{\prime}= S{}^{\prime}\cap Z(G)$.'' [The priming
here indicates commutator subgroup.]

In the focal subgroup paper [DGH3], p.\thinspace 496, there is a final remark about a
communication with Brauer who reports overlaps with results from his paper
on the characterization of characters \cite{Brauer}. Higman
states that those results follow from the focal subgroup theory.
David Gluck brought to our attention the footnotes in Section 5 of \cite{Brauer}
which indicate that Richard Brauer essentially agreed with Don's assertions. It seemed to Gluck that Brauer was impressed with Higman, and this may have influenced Don's success in getting a job at the University of Michigan (for Brauer had been a member of the  University of Michigan faculty 1948--51 before moving to Harvard).
There seem to be logical connections between the characterization of characters and the focal subgroup theorem which are not yet fully explored.

The focal subgroup definition bears strong
resemblance to the stable cocycle condition described by Cartan--Eilenberg \cite{ce}, which, as Alperin has suggested could be  viewed as a way to extend the focal subgroup theory to higher cohomology groups.
We have no evidence that Higman foresaw these connections.

\subsection{Representation theory}

Donald Higman wrote 8 papers on representation theory, over a period of seven
years 1954--1960. The most influential paper was the first [DGH5], in which he
introduced the notions of a relatively projective (or injective) module for
the group algebra of a finite group, and gave the famous criterion for
relative projectivity which bears his name, {\it Higman's
criterion}.
James A. (Sandy) Green, writes \cite{Green}, ``I am
very sorry to learn that Don has died. Mathematically I owe him for the idea
of relative projectivity of modules over finite groups. This was the starting
point of work I did on the `vertex' and `source' of an indecomposable
$kG$-module ($G$ finite group, $k$ field of characteristic $p>0$)." One of the first results of any kind on indecomposable modules was Don's
paper [DGH6] which showed that having a cyclic Sylow $p$-subgroup was a necessary
and sufficient condition for a group to have just finitely many isomorphism
classes of indecomposable modules over a given field of characteristic $p$.
The theory has come a
long way since then, and some of its main contributors have written in this
very volume. In particular, Michel Brou\'{e}  \cite{Broue} discusses modern
generalizations of Higman's criterion, very much in the spirit of Don's own
generalization [DGH7] beyond group algebras. Michel writes in the introduction of
\cite{Broue}, ``It must be noticed that ways to generalizations of Higman's criterion
had been opened half a century ago by Higman himself...".   Again in a general
algebraic context, rather than for group algebras alone, Don presented in [DGH9]
a theory, independently obtained by G. Hochschild \cite{hochschild}, of {\it relative homological
algebra.} 

     Don was keenly aware that not all interesting algebras were group algebras, and he was especially attracted to integral representations of orders over an integral domain O, in a semisimple or separable algebra (over the domain's quotient field). In the separable algebra case Don defined [DGH8] an ideal in O, now called the Higman ideal \cite{rhd}, pp. 253-258, \cite{CurtisReinerNew}, pp.603--609, as the (nonzero!) intersection of the annihilators of all first cohomology groups of bimodules for the order, an ideal which Don notes (see [DGH12]) to be contained in the intersection of the annihilators of all degree 1 Ext groups between lattices for the order.
If O is a Dedekind domain, the 1-cohomology and degree 1 Ext groups used are the usual ones, but require modifications otherwise [DGH12]. The intersection ideal defined with degree 1 Ext groups also makes sense in the somewhat weaker case of orders in semisimple, rather than separable, algebras.   (All semisimple algebras are separable over fields of characteristic 0.)
In the case of a finite group algebra over O, the principal ideal generated by the integer which is the group order is contained in all of these annihilator ideals. In this respect the Higman Ideal and its variations provide a substitute in the general case for the group order.  Higman used these ideals, and other variations defined for degree 1 Ext groups with a fixed lattice as one of the two module variables, to establish in [DGH8] and [DGH12] a series of results for general orders (in a separable algebra) over Dedekind domains that had been proved by Miranda \cite{miranda1,miranda2} for finite group algebras.
He and Jack McLaughlin also
studied orders in the function field case [DGH10].

Len Scott writes ``One of my favorites among Don's representation theory papers is [DGH11].
He proved that the issue of whether
or not two orders over a complete principal ideal domain were isomorphic could
be reduced to similar questions over a single finite length quotient of the
domain, and he showed that there is a bound on such a length which is
computable in terms of the Higman ideal. The proof is a convincing application of Don's homological ideas,
and the result itself is an important landmark in the theory of group ring
isomorphims, as well as of isomorphisms between more general orders.''

\subsection{Influence of the work of Donald Higman, by Charles Curtis}

Charles Curtis wrote to the editors as follows about Don's influence on his work, in particular regarding the widely-used reference 
\cite{CurtisReiner}:

``The book of Reiner and myself was intended to be an introduction to the
representation theory of finite groups and associative rings and algebras. One
objective was to bring readers to the point where they could approach the
papers in the area by Richard Brauer. While Brauer's work already contained
important contributions to ordinary and modular representation theory and
applications to the structure of finite groups, especially the finite simple
groups, we believed that this area offered rich possibilities for the future.

``One thing that Brauer, Nakayama, and Nesbitt had begun to develop was the
theory of nonsemisimple algebras, in particular group algebras of finite
groups over fields whose characteristic divided the order of the group. Two
chapters of our book were devoted to this subject. We tried to make full use
of the new concepts of projective and injective modules, introduced a short
time before by Cartan and Eilenberg. The adaptation of these ideas to the case
of group algebras was done by Gasch\"{u}tz in 1952, and extended to the
notions of relative projective and injective modules by Higman in 1954. The
main result obtained at this stage in the development was Higman's theorem,
published in 1954, that a group algebra of a finite group $G$ over a field of
characteristic $p$ has finite representation type (that is, has at most a finite
number of isomorphism classes of indecomposable modules) if and only if the
Sylow $p$-subgroups of $G$ are cyclic. This result has been influential in the
theory of integral representations of finite groups, in the representation
theory of Hecke algebras, and in the modular representation theory of finite
groups, up to the present time. Another important development was Green's
theory of vertices and sources of indecomposable modules, which remains a
central topic in the modular representation theory of finite groups, and is
based on Higman's theory of relative projective modules. One of the chapters
of our book was largely devoted to Higman's work in this area."

\section{Details of later work: algebraic combinatorics, groups and
geometries}

\subsection{ABA groups, rank 2 geometry characterization}

The ``ABA groups'' paper [DGH13] in 1961 with University of Michigan colleague Jack McLaughlin
arose from an effort to characterize low rank finite simple groups and their geometries.
Apart from the famous theorem that   a finite 2-transitive projective plane is Desarguesian, the paper contains
a number of fundamental results that continue to have application in group theory and combinatorics.
Shortly after it appeared, a result from the paper was used by Charles Curtis \cite{Curtis} to classify a family of
finite Chevalley groups, while nearly 30 years after its publication, another of its results  laid the basis for the classification \cite{Bueketal}
of flag-transitive linear spaces, proving that such spaces have a point primitive automorphism group. The paper
remains influential, being cited nearly 20 times in the last decade in publications on linear spaces and symmetric designs.

\subsection{Rank 3 permutation groups}
Don Higman's theory of rank 3 permutation groups was initiated in his
paper
[DGH15]
in 1964 in which he studied their parameters,
incidence matrices and character degrees. His approach was highly influential
on the work of other group theorists, not least because several sporadic
simple groups were discovered as rank 3 groups, in particular the Higman--Sims
simple group discussed below. Also in his 1964 paper Higman gave a pair of design constructions which are now standard constructions for symmetric designs.
His other papers on this theme focused on characterisations proving,
for example, that all rank 3 affine planes are
translation planes.
It is notable that Higman chose rank 3 groups
as the topic for his invited lecture
[DGH28]
at the International Congress of
Mathematicians in Nice in 1970.

\subsection{The Higman--Sims graph and the Higman--Sims sporadic simple group of
order $44,352,000=2^{9}3^{2}5^{3} 7{\cdot}11$}

The discovery by Donald Higman and Charles Sims of the Higman--Sims simple
group is legendary. According to Sims' account in \cite{Hiss}, Don
and he had just heard Marshall Hall's description, at the 1967 ``Computational
Problems in Abstract Algebra'' conference in Oxford, of the construction of
the Hall--Janko sporadic simple group as a rank 3 permutation group on 100
points. Higman and Sims examined closely a list of possible rank 3 parameters
that Higman had compiled using his rank 3 theory, and quickly focused their
attention on a possible rank 3 group with point stabiliser the Mathieu group
$M_{22}$ or its automorphism group $M_{22}.2$. In discussion during an
intermission in the formal conference dinner on the last day of the
conference, September 2, 1967, Higman and Sims realised that these groups had
natural actions on both 22 points and 77 points. They continued their work on
these actions after the dinner finished, and eventually realised that they
needed to construct a valency 22 graph on 100 vertices, now called the
Higman--Sims graph. In the early hours of Sunday, September 3, 1967, using the
uniqueness of the associated Witt design, they proved that the automorphism
group of their graph was vertex-transitive, and was a new simple group.

We re-publish with permission the 2002 account by Charles Sims \cite{Hiss} of this famous story.

``Prior to this conference, Don had been investigating rank 3 groups. He had
discovered a number of conditions that the parameters [DGH57]
of such a group have to
satisfy and had used a computer to generate a list of parameters that
satisfied all of his conditions. I was familiar with Don's work. [\dots]

``At the Oxford conference, Marshall Hall announced the construction of Janko's
second group.\footnote{This group is more correctly  called the Hall--Janko
group since both Hall and Janko discovered it independently; see \cite{GreviewR,Griess12}.} There is a long paper by Marshall in the
conference proceedings that includes, among other things, a description of how
this group was constructed. The group was given as a rank 3 group of degree
100 with subdegrees 1, 36 and 63.

``After Marshall's talk, Don and I discussed whether there might be other rank 3
groups of degree 100. If it were not the case that we use the decimal system
and that $100=10^{2}$, I am not sure we would have asked this question.

``Don consulted his list of rank 3 parameters and found
that the subdegrees 1, 18, 81 appeared on the list of degree 100. It did not
take us long to realize that the wreath product of $S_{10}$ and $Z_{2}$ has a
rank 3 representation (on the Cartesian product of two copies of $\{1,\cdots,
10\}$) of degree 100 with these subdegrees.

``Encouraged, we looked at the list again. There was one other set of subdegrees
for a possible rank 3 group of degree 100, namely 1, 22, 77. The number 22
certainly suggests that the stabilizer of a point should be $M_{22}$ or its
automorphism group. ($S_{22}$ and $A_{22}$ don't have representations of
degree 77.) We did not immediately see how to construct such a group and
agreed to continue our discussion later.

``The conference dinner was held on the last day of the conference, Saturday,
September 2. This dinner was quite formal. After the main part of the meal,
the participants were asked to leave the hall while the staff cleared the
tables and prepared for dessert and coffee. As Don and I walked around the
courtyard of the college in which the dinner was being held, we again talked
about 100=1+22+77. The first question to answer was whether $M_{22}$ has a
transitive representation of degree 77. Both Don and I were familiar with
combinatorial designs and knew that $M_{22}$ acts on $\mathfrak{S}(3,6,22)$.
As we walked, we computed the number of blocks in this design. When we got the
answer 77, we were sure we were on to something.

``At this point, it was time to go back in for dessert. After the dinner was
finished, we went to Don's room [\dots] to continue working. There were some
false starts, but eventually we realized that we needed to construct a graph
of degree or valence 22 with 100 vertices consisting of a point *, the 22
points of an $\mathfrak{S}(3,6,22)$ and the 77 blocks of that $\mathfrak{S}%
(3,6,22)$. The point * would be connected in the graph to the 22 points. One
of the 22 points $p$ would be connected to * and to the 21 blocks containing
$p$. A block would be connected to the 6 points in the block and to 16 other
blocks. We had to do some computations, but it was not hard to show that a
block is disjoint from exactly 16 other blocks. Thus in the graph, two blocks
should be connected if they are disjoint.

``At this point, we had the graph and we knew that $Aut(M_{22})$ as the
stabilizer of * in the automorphism group of the graph, but we did not know
that the automorphism group had any elements that moved *. To get additional
elements, we used the fact that Witt had proved that $\mathfrak{S}(3,6,22)$
was unique.

``I don't think we finished until the early morning hours of Sunday, September
3, 1967. There was one uncertainty. We knew that we had a new simple group,
but we did not know whether the stabilizer of a point in the simple group was
$M_{22}$ or $Aut(M_{22})$.
[\dots] ''

\subsubsection{Mesner thesis and Higman Sims graph}

While researching this biography, the editors learned that the 1956 doctoral thesis of 
Dale Mesner
\cite{mesner}
 contains
a construction of the Higman--Sims graph; see also \cite{jj}.   
Higman and Sims were unaware of the
Mesner thesis.  
Jon Hall kindly gave us the following summary of
its contents.

``This long thesis (291 pages) explores several topics,
including integrality conditions for strongly regular
graphs (association schemes with two classes)
related to Latin squares. One case of particular interest
is that of graphs with the parameters of the Higman-Sims graph.
Mesner notes that the existence of such a graph
is equivalent, via a graph construction essentially
the same as that of Higman and Sims, to the
existence of a balanced incomplete block design
with the parameters of the $S(3,6,22)$ design
(including having 16 blocks disjoint from each
block). 
He proceeds to construct such a design.  
At
a certain point in this construction he is faced 
with four choices, among which he cannnot
distinguish. He makes one choice and completes
his construction, saying that the other three cases
give similar results. At the end he then has
four constructions of graphs on 100 points having
the same parameters as the Higman-Sims graph.
Of course, what he has is four slightly different
constructions of the same graph, but he does
not know or investigate that. On page 147 he 
finishes this part of the thesis by saying, 
`This completes the construction of scheme \#94 [the graph]. 
There are at most four solutions to the association scheme,
corresponding to the four choices for the structure of the
blocks of [the design]. It is not known whether any of the
four solutions are equivalent under some permutation of
treatments.' 
There is no further determination of isomorphisms or the automorphism group, 
nor mention of the groups of Mathieu or the work of Witt on Steiner
systems."


\subsection{Strongly regular graphs and association schemes}

Don's original motivation for studying the combinatorics underlying
permutation groups was purely group theoretic. For example,
in his  magnificent ``intersection matrices'' paper
[DGH18]
of 1967,  Don
proved his famous primitivity criterion for
finite transitive permutation groups mentioned in the introduction to this article.
In that same paper he introduced (what are now called)
{\it intersection arrays} of finite distance regular graphs from a group
theoretic perspective, and began the study of  permutation groups of maximal
diameter.
Gradually it was realised that rank 3 permutation groups
correspond to strongly regular graphs, that Don's intersection
matrices correspond to association schemes,
and  permutation groups of maximal
diameter correspond to distance transitive graphs,
a class of graphs studied intensively over the succeeding decades.
The concept of a strongly regular graph
already existed in the combinatorial literature,
and association schemes (that is, symmetric association schemes)
were known in
statistics and experimental design theory, in particular in the work of R. C.
Bose and his school.

It seems that these connections were not known to group theorists
at the time of Don's paper
[DGH18].
Several years later in 1970,
in an important paper
[DGH26]
with Hestenes, Don
explored the connection between rank 3 groups and strongly
regular graphs, introducing and studying the adjacency algebras of such
graphs and giving a combinatorial form of primitivity.
The theories of strongly regular graphs and
association schemes were given real mathematical depth by this
fruitful encounter with group theory.

Moreover it was natural for Don to study these objects also
from the combinatorial viewpoint: in his 1970 paper
[DGH26],
he introduced the ``$4$-vertex condition'' on a graph, a combinatorial
alternative to the group theoretic rank $3$ condition. He used this to obtain
a new proof of a result of Seidel determining the strongly regular graphs with
smallest eigenvalue $-2$.

\subsection{Coherent Configurations}

One of Don's most notable contributions in the area of combinatorics is his
introduction of the concept of a \emph{coherent configuration} in his 1970
paper
[DGH23].
A
coherent configuration is an
axiomatization, in a purely combinatorial setting, of the structure
of an arbitrary  permutation group. It generalises the notion of an
association scheme, which may be regarded as a combinatorial
axiomatization for a transitive
permutation group. Don arrived at this general concept as soon as he
was aware of the combinatorial interpretation of his early work on
permutation groups. He
developed the theory of coherent configurations in full
generality: first in~
[DGH23] 
and in more polished form in
other papers, most notably 
[DGH35], [DGH37]
and his paper [DGH41]
on coherent algebras in 1987. The main
content of these papers is a study of the algebraic conditions arising
from association schemes and coherent configurations. Don studied the
adjacency algebras (centralizer
algebras) of these combinatorial structures, obtaining many important properties
such as {\it orthogonality
relations, Schur relations}, and {\it Krein conditions}. \\

A major reason why Don studied the representation theory of association schemes and
coherent configurations was to
apply the algebraic conditions he obtained to the study of many geometric objects.
In his 1974 paper
[DGH32]
he illustrated how the theory of
coherent configurations could be applied by  studying
{\it generalized polygons}, and in particular giving alternative proofs of the
famous theorem of W. Feit and G.~Higman, and of the  Krein condition,
discussed below. These results, in turn,  were applied in
[DGH32] 
to prove the following theorem: if a generalized quadrangle or octagon
has $s+1$ points on each line and $t+1$ lines through each point, with
$t>1$, then $s\leq t^2.$ Don also applied his theory more generally, among other things,
to study {\it partial
geometries}, a concept related to that of a strongly
regular graph.\\

Don Higman's representation theory of coherent
configurations is thorough and useful, and unequalled by any other author.
On the other hand, he did not treat
many concrete examples in his papers, leaving this to others such as \cite{bannaiito}
and \cite{bcn} (for symmetric or commutative association schemes).\\

The non-negativity of the so-called Krein parameters
is very important in the theory of association schemes, and in particular in the
Delsarte theory of codes and designs in association schemes \cite{bannaiito,bcn,Delsarte}.
A version of this condition was first noticed in 1973 in the context of permutation groups by Leonard
Scott and his colleague Charles Dunkl, translating from earlier work of Krein, see~\cite{Scott} and
\cite{Scott2}.
Using work of Schur, dating even earlier than Krein's work, Don (and independently
N. Biggs~\cite{Biggs}) proved the Krein condition in full
generality in the context of association schemes.
Don's new approach to the development and exposition of
the Krein condition theory has been very influential \cite{Coolsaet,Gavrilyuk,Hiraki,Hobart}.   
Don often lectured on this to new audiences.

\subsection{Geometric applications of coherent configurations}

Donald Higman's  research in 1980's and 1990's focused mainly
on the application of the theories of association schemes and coherent configurations to geometry.
Many of his later works may be characterised as follows. Using the algebraic properties of
adjacency algebras of association schemes and coherent
configurations, he determined
their feasible parameter sets under various conditions. He then characterized the known
examples, and for many of the other small parameter sets, he proved the nonexistence
of examples. He was interested in a variety of geometric structures, especially those
admitting interesting group actions.
In particular, he studied special kinds of buildings (in the sense of Jacques Tits)
as imprimitive association schemes. We discuss below other notable geometric structures
he studied.

\medskip\noindent
A transitive
permutation representation can be regarded as the induced representation of the identity 
representation of a subgoup $H$ of $G$. A similar theory is available for
the induced representation of a linear representation of a subgoup $H$ of $G$.
Don's general theory of weights on coherent configurations (and association schemes) is
a combinatorial analogue of this. Using this theory Don gave combinatorial proofs
in 
[DGH33], [DGH37]
of several group theoretic
results due to Frame, Wielandt, Curtis--Fossum, and Keller (see \cite{CurtisFossum,Keller,Wielandt}).
An especially nice result of this type
is his calculation in 
[DGH36] of the
degrees of the irreducible representations appearing in the induced representation
of the alternating character of the subgroup $2\cdot {}^2E_6(2)\cdot 2$
of the ``Baby Monster'' sporadic simple group.

In [DGH41], 
 Don introduced several geometric structures associated with small rank
coherent configurations and proposed that they be studied systematically. Typical examples
are those with the following four `Higman parameters':
\begin{center}
 (i)\quad $\binom{2 2}{\phantom{2} 2}$,\quad (ii)\quad  $\binom{2 2}{\phantom{2} 3}$, \quad (iii)\quad
  $\binom{3 2}{\phantom{2} 3}$,\quad (iv)\quad  $\binom{3 3}{\phantom{2} 3}$
\end{center}
where the corresponding geometric structures are (i) symmetric designs,
(ii) quasi-symmetric designs,
(iii) strongly regular designs of the first kind, and
(iv) strongly regular designs of the second kind.
Much work had already been done in cases (i) and (ii),
and Don studied cases (iii) and (iv)  in  
[DGH42, DGH48]. In the meantime Hobart \cite{Hobart}
studied another case  $\binom{2 2}{\phantom{2} 4}$. Again, the techniques are to study algebraic
properties, and obtain a list of feasible parameters of small sizes.

Motivated by examples associated with classical
triality and some other related group theoretic examples,
Higman [DGH47]
studied imprimitive rank $5$ permutation groups.
He divided them into
three cases according to the rank and corank of a parabolic subgroup, and
studied algebraic conditions of the associated imprimitive association schemes.
His student Yaotsu Chang had earlier dealt with the imprimitive rank $4$
case in his Ph.\thinspace  D. thesis.  In an unpublished
preprint [DGH55] Don extended his research to
include uniform $(t,p,m)$ schemes, in particular,
such structures afforded by a transitive action of a group of the form
$G \cdot S_{t+1}$ with $G$ simple and $S_{t+1}$ acting faithfully
on $G$. He focussed on the cases $t=2, 3$, and eliminated those
with small sizes by studying algebraic conditions of their associated coherent
configurations. Again, group theoretic assumptions seem necessary to
obtain complete classifications.

In  [DGH45] and an  unpublished preprint [DGH52],
Don studied regular $t$-graphs, as a generalization of
regular 2-graphs in the sense of Seidel. In particular, a regular 3-graph is naturally associated
with a coherent configuration of rank 6. By studying their algebraic
properties, he obtained a list of feasible small parameters for such
coherent configurations, and then studied whether such configurations exist.
However, a complete classification of regular $t$-graphs for large $t$,
or even for $t=3$, seems infeasible.
The Ph.\thinspace  D. thesis of Don's student Sankey in 1992 discusses regular
weights on strongly regular graphs, a different generalisation of $2$-graphs,
but still in the spirit of [DGH45].

In concluding this subsection, we wrap up the correspondence between groups
and combinatorics in a small dictionary as follows:
\begin{center}
rank 3 group $\leftrightarrow$ strongly regular graph

transitive permutation group $\leftrightarrow$ association scheme

arbitrary permutation group $\leftrightarrow$ coherent configuration

distance-transitive group (graph) ( = permutation group of maximal
diameter) $\leftrightarrow$ distance-regular graph.
\end{center}

In passing, we would like to mention the following two streams of similar ideas
to study combinatorial objects from a group theoretical view point.

(1) Essentially the same concept as
coherent configuration was studied in Russia in the name of ``cellular rings",
independently and simultaneously, as there was not much scientific
communications between east and west at that time. (See, for example, \cite{weisfeiler}.)

(2) The theory of Schur rings developed by Wielandt \cite{Wielandt} and then further
by Tamaschke  \cite{tamaschke1}  also treats a connection between
groups and combinatorics. (Here, note that the concept of Schur ring is
something between permutation groups and association schemes.)

\subsection{Connections with statistics, topology, mathematical physics}

There are several connections which the work of Donald Higman makes outside algebra.
Association schemes occur in statistics (this article mentions the work of Bose and Bailey).
The thesis of Mesner \cite{mesner} on Latin squares turned up a special set of graphs on 100 points, one of which is the Higman--Sims graph.
Fran\c cois Jaeger constructed the spin model (in the sense of V. H. F. Jones)
on the Higman--Sims graph, and it plays a very important role in the theory of
spin models (more broadly in topology and mathematical physics.)
See \cite{delaHarpe},  \cite{Jaeger}.

\section{Younger mathematicians and Donald Higman}

Don Higman's international engagement in the group theory and
combinatorics communities made a strong impression on younger
mathematicians. He made regular visits to England, Holland,
Belgium, Germany, Australia, Italy and Japan and organized
meetings in Oberwolfach starting with the meeting  ``Die
Geometrie der Gruppen und die Gruppen der Geometrie under
besondere Ber\"{u}cksichtigung endliche Strukturen", 18--23
May, 1964.

\subsection{Michael Aschbacher}

``I went to graduate school to become a combinatorist, and wrote my
thesis in design theory. My approach was to study designs from the
point of view of their automorphism groups, so I became interested
in finite permutation groups, to the point that in my last year in
graduate school I decided to work in finite group theory rather
than combinatorics.

``In one part of my thesis, I studied certain designs whose
automorphism group was rank 3 on points of the design. Thus I was
led to read Donald Higman's fundamental papers on rank 3
permutation groups, and rank 3 groups have remained one of my
interests ever since.

``In particular in one of my early papers \cite{Asch}, I showed that there is no
Moore graph on 3250 vertices with a rank 3 group of automorphisms.
This was motivated by Higman's paper 
[DGH15],
where in section 6
he shows that in a rank 3 group with $\lambda=0$ and $\mu=1$, $k$
is 2, 3, 7, or 57, and he determines the possible groups in all
but the last case, which he leaves open.
In the graph theoretic literature, the corresponding strongly
regular graphs are called \textit{Moore graphs}, and of course a
Moore graph of valence 57 is on 3250 vertices.

``A year later, reading Bernd Fischer's wonderful preprint on groups
generated by 3-transpositions, I decided to attempt to extend
Fischer's point of view to what I called `odd transpositions'.
The base case in this analysis involves almost simple groups $G$,
and after handling some small degenerate cases, an extension of
Fischer's fundamental lemma shows that $G$ acts as a rank 3 group
on a certain system of imprimitivity for the action of $G$ on its
odd transpositions.
Thus, once again, Higman's theory of rank 3 groups became a
central tool in my early work.

``I probably first met Donald Higman at a meeting in Gainesville
Florida in 1972, organized by Ernie Shult. However my strongest
early memories of Higman are from a two week meeting in Japan in
1974, organized by Michio Suzuki. This was only my second trip
abroad, and one of my favorites. A number of wives also took part,
including my wife Pam. Both Donald and Graham Higman attended the
meeting, and I recall Pam referring to them as ``the Higmen". As I
recall the two were connected in our minds, not just because of
their names, but because both had extremely luxuriant eyebrows.
....

``He was an excellent mathematician who will be missed.''

\subsection{R. A. Bailey}

Rosemary Bailey relates three small anecdotes about Don, whom she, like many of her fellow research students in Oxford, called DGH.

``I started as a Ph.\thinspace  D. student at Oxford in October 1969, under the
supervision of Graham Higman. In summer 1970 I attended the
International Congress of Mathematicians in Nice. Of course, much
of it was above my head at that stage. The talk that most
impressed me was the one by DGH. I was even brave enough to go up
and talk to him about it afterwards, and I remember him being very
kind.

``Some time in 1971, DGH visited Oxford for a few months. He gave
a course of lectures called `Combinatorial Considerations about
Permutation Groups'. At the time, not only was I working on
permutation groups but Graham Higman had given us a series of
advanced classes on strongly regular graphs, so we were well
prepared. I was absolutely fascinated by DGH's clear lectures.
Well, I went away and became a statistician, but something must
have stuck in my mind, because when I came back into association
schemes via the original Bose--Nair approach I was able to put it
all together, and this eventually led to my 2004 book \cite{Bailey} on
Association Schemes. I mentioned these lectures of DGH's in the
Acknowledgements page to that book.

 ``Most combinatorialists know what a $t$-design is, or at least
what a 2-design is. Dan Hughes introduced this terminology into
the literature. However, he says that it was DGH who suggested it
to him: see pages 344-345 of \cite{Bailey}.  ''

\subsection{Eiichi Bannai}

``I knew Donald G. Higman for many years, and I was mathematically
very much influenced by him. This was at first through Suzuki,
whom I viewed as a kind of mathematical father. In an article
written in Japanese, Suzuki wrote as follows:  `When I (Suzuki)
went to Illinois in 1952 Donald Higman was there and just
completed his Ph.\thinspace  D. thesis on focal subgroups. Partly because I
knew the name of Higman already, by his paper on homomorphism
correspondence of subgroup lattices, soon we became very good
friends.' I also remember that Suzuki described Higman as a
mathematician doing very fundamental and essential work with his
own viewpoint.

``Sometime, in the 1960s, Higman visited Japan and stayed with the
family of Suzuki (who had settled at Illinois, but regularly
returned in the summers--Suzuki and Higman were sufficiently close
that I came to regard Higman as a kind of uncle). \ I found an
official 3 page note published in Sugaku (the official journal of
Mathematical Society of Japan) which summarizes Higman's talk on
June 3 in 1966 at the University of Tokyo. This note was prepared
by Koichiro Harada, then a young graduate student. The title of
his lecture was: Remarks on finite permutation groups. According
to the note, Higman talked on the discovery of Janko's new simple
group, and discussed D. Livingstone's new construction of it as
the rank 5 permutation group of degree 266 with subdegrees $1,\
11,\ 110,\ 12,\ 132.$ Then he (Higman) discussed several attempts
of his own and of McLaughlin to try to construct (new) rank 3
permutation groups. In particular, he mentioned his result that if
there exists a rank 3 permutation group with subdegrees $1,\ k,\
k(k-1),$ then $k=2,\ 3,\ 7,$ and possibly 57. (At that time, it
seems that he was not aware of the connection with strongly
regular graphs and the work of Hoffman-Singleton.) In addition, he
discussed some of the work of Tutte, Wielandt, Sims, etc.
Interestingly enough, the Higman--Sims group HS was not mentioned at
all, as it was before its discovery.

``The content of this note tells us that  Higman had several years
of serious preparations before his discovery of  the group HS. The
discovery was clearly not just a matter of luck. On the other
hand, there was also some element of good fortune. See the account
of Sims in section 5.3.

\ \ \ \ \ ``Higman also visited Japan for the Sapporo Conference in
1974 (where I first met him in person). Since then I had many
chances to meet him in the USA and in Europe, and he was very
generous and helpful to others. He contributed to a good research
environment in the field and greatly helped other researchers in
many ways, for example, by organizing several Oberwolfach
conferences, etc.. When Suzuki visited Japan during summers, he
always gave a series of lectures on very hot new developments of
group theory. ... I was able to attend Suzuki's lectures for the
first time in 1968. In his lectures, Suzuki covered many topics,
but the main topic was the discoveries of new finite simple
groups, including Hall-Janko, Higman--Sims, and several others. ...
I learned the charm of rank 3 permutation groups from Suzuki's
lectures. At that time, I was studying representation theory of
finite groups (through the book of Curtis-Reiner and also from
Iwahori, my advisor) and also finite permutation groups, in
particular multiply transitive permutation groups (through the
book of Wielandt \cite{Wielandt}). After Suzuki's lectures, I started to study
rank 3 permutation groups, along the line of Donald Higman. Also,
I was interested in the combinatorial aspect of finite permutation
groups, such as strongly regular graphs. The connection between
strongly regular graphs and rank 3 permutation groups was well
known by that time. I read many of Higman's papers on rank 3
permutation groups, those related to classical geometries, and the
one on intersection matrices for finite permutation groups, etc.
Although Higman's earlier papers were in the context of finite
permutation groups, the combinatorial contents in it was very
evident.\newline \ \ \ \ \ In early 1970's, I tried to follow
Higman's research direction in various ways. Let me comment on
some of these directions of Higman:\newline

``(i) (Small rank subgroups of classical groups.) The first thing I
tried in my research was to find  (and classify) small rank (in
particular rank 2, i.e., 2-transitive)  subgroups of classical
groups, or of simple groups of Lie type. During that study, I
determined rank 2 subgroups $H$ of $PSL(n,q)$, i.e., 2-transitive
permutation representations of $PSL(n,q)$, by determining the
candidates of possible irreducible characters
$\chi$ which might possibly appear in
$1_{H}^{G}$, I used Green's theory of the character theory of $GL(n,q)$.  Later, I  realized that this approach is
considerably simplified by using the work of Higman  [DGH14] on flag
transitive subgroups, since it holds that $\chi$ must appear in
$1_{B}^{G},$ if $H$ is not flag transitive. (Here $B$ is the Borel
subgroup, i.e., the  upper triangular subgroup.)  Also, I was very
much interested in the paper [DGH13] which, before Tits general
theorem on spherical buildings,  characterized the group
$PSL(3,q)$ by, essentially giving the classification of $A_{2}$
type buildings (in the context of groups). This was used in many
situations later by many authors.

``(ii) (Rank 3 transitive extensions.) Another direction of my
resarch was to try to find new simple groups, starting with a
known group $H$ to find $G$ which contains $H$ as a rank 3 (or
small rank) subgroup. To be more precise, starting from a
permutation group $H$, to find a transitive group $G$ whose
stabilizer of a point is isomorphic to $H$ is called a transitive
extension problem. The discovery of the Higman--Sims was one such
beautiful example. Various sporadic groups, such as Hall-Janko,
McLaughlin, Suzuki, Rudvalis, and the Fischer groups, were
discovered this way. I think I started this project too late, as
all were already discovered, as it turned out, a consequence of
the classification of finite simple groups.. \newline\ \ \ \ \
(iii) (Intersection matrices for finite permutation groups.) The
paper I read most carefully was the paper [DGH18] with the title:
Intersection matrices for finite permutation groups. That paper
gave a systematic study of finite permutation groups with some
good property. The paper also gave essentially a new proof of
Feit--Higman' theorem on generalized polygons (in the context of
groups). This new proof, which was to look at $1_{P}^{G}$ rather
than $1_{B}^{G}$ in the original proof, was very transparent for
me, and I was very much impressed. (Here, $P$ is a maximal
parabolic subgroup.) \ I noticed that many of the algebraic
properties of association schemes were also in this paper.
Actually, I can say that I first learned association schemes
mainly from this paper of Higman. In that paper he used the
notation for intersection numbers as $A_{i}A_{j}=\sum
p_{j,k}^{(i)}A_{k}$ [and actually, such an explict product is not
written down in the paper] which was somewhat different from the
traditional one $A_{i}A_{j}=\sum p_{i,j}^{(k)}A_{k}$ used in the
theory of association schemes [apparently unknown to Higman at the
time]. I was so accustomed to using Higman's notation that I first
insisted using his notation when I later wrote a book \cite{bannaiito} with
Tatsuro Ito, though we later returned to the traditional notation
(as did Higman). Anyway, the influence of Higman's intersection
matrix paper on me was  very great.  \newline

``There are many other areas of influence on my research, and I have
described them in contributing to sections 5.4 (strongly regular
graphs and association schemes), and coherent configurations
(section 5.5). Let me add here that Higman's development of the
general theory of coherent configurations in his papers
[DGH35], [DGH37] is
very thorough and beautiful. Though not many concrete examples are
treated, his papers are also useful as well as complete. There is
nothing comparable to the work of Higman in the study of coherent
configurations.

``Later, Higman did study many concrete examples of coherent
configurations with special types, in papers [DGH41], [DGH42], [DGH45], [DGH47], [DGH52], [DGH53], [DGH54], [DGH55], etc. They are
related to some kinds of designs and finite geometries. While
these studies may not be breakthroughs themselves, they indicate
many possible future research directions. Some of the papers in
this special volume will reflect the Higman's strong influence
along this line of study. Coherent configurations will be studied
more in the future, as they are a very natural concept. Yet, it is
difficult to know what is the most important general research
direction for their study (beyond association schemes). I
personally feel that Higman's research in late years was to try to
find out some new possibilities, by trying several directions.

``In conclusion, let me say Higman's work since 1960 has been
extremely influential on those of us working on groups and
combinatorics. \ Higman opened a new path which led to the
development of so called algebraic combinatorics, and we regard
him as one of the founders of this approach. Higman continued his
research in more combinatorial context, by giving solid
foundations of the theory of coherent configurations, and also
studied geometric (design theoretical ) aspects of the theory of
coherent configurations. He succeeded in combining group theory
and combinatorics very nicely, and he will be remembered by the
mathematical community for his fundamental contributions to these
areas. Also he will be remembered for his kind and generous
personality by those of us who know him personally.''

\subsection{Peter J. Cameron}

``I arrived in Oxford to do my D.Phil. (Oxford term
for Ph.\thinspace D.) in 1968, and worked under the supervision of Peter
Neumann. I worked on permutation groups. At the time, the most
exciting development was the introduction of combinatorial
techniques into the study of permutation groups, due to Charles
Sims and Donald Higman. Sims's methods were graph-theoretic, based
on the work of Tutte; Higman's were more representation-theoretic,
and grew from the work of Wielandt. All my subsequent work has
grown out of my initial exposure to this material.

``It was very important to me that Don Higman spent part of the year
1970--1971 in Oxford, for several reasons. First, he gave a course
of lectures. I was one of two students (the other was Susannah
Howard) responsible for producing printed notes after each week's
lectures, and compiling them into a volume of lecture notes
published by the Mathematical Institute, Oxford, under the title
"Combinatorial Considerations about Permutation Groups". These
notes contained a number of ideas which only appeared later in
conventional publications; I was in the very fortunate position of
having a preview of the development of the theory of coherent
configurations.

\dots

``Fourth, and most important for me. Don arranged for me to spend a
semester at the University of Michigan as a visiting assistant
professor. This was a very productive time for me. I was allowed
an hour a week in which I could lecture on anything I liked, to
anyone who wanted to come along. I did a vast amount of
mathematics there, and made many good friends.''

\subsection{Robert L. Griess, Jr. }

``When I joined the University of Michigan Mathematics Department faculty in September, 1971, I met Don Higman for the first time.  I knew his name because the simple group constructed by Don Higman and Charles Sims had been well-studied by members of the group theory community for a few years.  About two or three years before at the University of Chicago, where I was a graduate student, I met Doris and Jack McLaughlin during Jack's sabbatical in Chicago, where he presented his new construction of the McLaughlin sporadic simple group.   The rank 3 theory of Higman was prominently featured in McLaughlin's lectures.  Roger Lyndon's achievements in combinatorial group theory, homological algebra and logic were well-known.  For these reasons, I felt attraction to the University of Michigan.

``When I met Don, his era of algebraic combinatorics had begun.
An important undercurrent was representation theory of groups and algebras. His productivity during the seventies and the attention it drew from experts worldwide were quite impressive.

``Don and I together attended hundreds of classes and seminars.  While we never collaborated on any research, we conversed for decades about algebra, especially representation theory, and  aspects of finite groups and combinatorics.  During this time,  I had many opportunities to see how his mind operated.

``A lot of Don's research involved long-term  methodical work.  He collected data which he studied as he developed his theories.  He was repeatedly attracted to the idea of elegant, simple explanations and to finding axiom systems.  His publications may not have clearly indicated the lengthy
efforts he made.

``About a year after arriving, I heard Don talk about his short result on transitive extensions of permutation groups [DGH30].
It was presented in Don's characteristic modest and humorous style.   It struck my young mind as quite original and fresh, a real piece of mathematical wit.   While I was a graduate student, I had gotten the impression that his work in finite group theory was fundamental.  Finite  group theorists  valued  the focal subgroup theorem, basic results on indecomposable modules, the idea of relative projectivity, and work on rank 3 permutation groups.

``Don regularly had mathematical visitors, including many who stayed at his home.  Betty and Don were generous hosts on many occasions and made me,  in particular,  feel very welcome during dinners and parties.
I appreciated the chance to get acquainted with so many mathematicians from Don's world.

``Betty and Don took friendly interest in my personal life and Don was always keen to hear about my mathematical activities.   He made helpful suggestions and on many occasions he indicated publications or preprints I did not know about which might lead to interesting connections.  I learned new viewpoints and details about the genesis of ideas in groups and geometries.

``Donald Higman was a model of perseverance, clarity and
elegance in mathematics.  His death in 2006 was a complete surprise.  Fortunately, I had dedicated an article to him a few years before \cite{bwy,bwycorr}.  When I told him, he was happy and thanked me quietly.  I shall never forget him.''

\subsection{Willem Haemers}

``When I met Don (probably in 1976), I was a Ph.\thinspace  D. student at the
Technical University in Eindhoven (supervised by Jaap Seidel).
Then Don gave a beautiful course on  `classical groups' ,
and I had the pleasure of participating in several conversations
with Don and Jaap.

``Already in those days Don had the habit of writing unfinished
manuscripts.
He often gave such manuscripts away, hoping that something would
come out of it. Sometimes this worked. But probably many such
manuscripts remained unpublished, It is how my two [joint] papers
with Don came to existence. \dots He always showed great interest
in me, in my mathematics and private life. I liked him very much. \footnote{The unpublished manuscripts and other material known to the editors are listed in this article.}''

\subsection{Cheryl E. Praeger}

``I was in Oxford for the course given by DGH and alluded to in the
writings of Ro Bailey and Peter Cameron. The course made a big
impression on me, and the papers by DGH on rank 3 groups and
intersection matrices were ones that I dissected more than read.
They seemed like magic - and those papers, together with papers by
Peter Cameron and Charles Sims that built on the notion of orbitals for transitive
permutation groups, and Norman Biggs's book on Algebraic Graph
Theory, were formative in my thinking and development in work on
automorphism groups of graphs.

``\dots I visited U Michigan for exactly one day in 1974. I remember
meeting DGH then (and also meeting Bob [Griess] for the first
time). DGH asked me to send him copies of my preprints. I cannot
explain how overwhelmed I felt by his asking me to do this. I
really felt like I was unworthy, that I should presume to send
this great man something that I had written. I know that sounds
bizarre, but that's how it felt to me in 1974.''

\subsection{C. C. Sims}

Charles Sims's first course in abstract algebra, when he was an undergraduate at the University of Michigan, had been taught by Don Higman.  Thus, Don was a teacher who became a colleague.  Charles writes as follows about his first consultation with Don.

``I think the first research interaction with Don occurred shortly
after I got to Rutgers [about 1966]. I read one of Don's papers
on rank 3 groups and as I remember it, Don asked about a group of
degree 50 with subdegrees 1, 7, and 42. I wrote to him to point
out that two graph theorists, Hoffman and Singleton, had proved in
1960 the uniqueness of the associated strongly regular graph
without the assumption that the graph had an automorphism group
that was edge transitive. This started the contacts that
ultimately [led] to the discovery of Higman--Sims.''

\subsection{Stephen D. Smith}

Steve Smith writes: 
``In the mid 70s, conversations with Don and Bruce Cooperstein
   (I think at Ann Arbor) showed me the relevance of Dynkin
   diagrams to the geometry of the Lie type groups (this was
   before I understood much about Tits buildings).  That in
   turn was the background for my input into the joint work
   with Mark Ronan on the 2-local geometries for sporadic groups.  \cite{rs} 
   (And hence for my later long-term interest in geometries for groups.)''

\subsection{John G. Thompson}

In the 90s, Don Higman told Bob Griess that young John Thompson
spoke to Don about some of his insecurities during his struggle to
prove the Frobenius conjecture. Thompson's eventual solution
resulted in his celebrated 1959 thesis and energizing of the
classification work on finite simple groups. Don also reported
that at one Oberwolfach meeting, he observed John Thompson pacing
for hours, pondering a step in the classification of finite simple
groups. After the pacing, Thompson was pleased to tell Don that he saw what to do.
The date of this event is unknown. However it probably refers to work on the
N-group Theorem, which places it in the mid to late 1960s.
Don's role as confidant confirms the general impression he
gave as a concerned and friendly authority figure in the groups
and geometry community.  

John Thompson says that this paragraph ``touches the right bases'' and 
adds that  ``One of my abiding memories of Don is that he came up with a simple group on  
100 letters\footnote{Don Higman and Charles Sims are co-discoverers.}
which was not the one found by Hall and Janko.  This really made me scratch my head and
wonder if the list of sporadic simple groups would terminate.  The equations
100 = 1 + 36 + 63 = 1 + 22 + 77 were really scary.''

\subsubsection{About the The Odd Order Theorem}

This is the famous and important theorem by Walter Feit and John Thompson
\cite{FeitThompson}, which says that finite groups of odd order are solvable.
It is technically very difficult and takes up
an entire issue of the Pacific Journal of Mathematics, 1963. Don
Higman told Bob Griess that he (DH) had been one of a team of
referees for that paper of 255 pages.   Don said that his
responsibility was mainly the generators and relations section at
the end.

Bob Griess acquired two copies of that issue from Don's
collection. One includes a handwritten note in Walter Feit's handwriting, 
on Yale University stationery.  It says only ``You are incorrigible. " It
is not addressed to anyone.  
It is in Don's possession, so was probably inserted in
an offprint which Walter gave Don.

\subsection{Students and associates}

The fifteen students of Don Higman are listed later in this
article. Many got jobs at research universities and published
actively. During the years of Higman and McLaughlin at the
Unversity of Michigan, these two often had many students and
there was a lot of interaction between them. The students of Roger
Lyndon, and Bob Griess, while in group theory and interactive with
the Higman and McLaughlin circle, were less numerous and
mathematically a bit distant.

Three of Don's students contributed articles to this memorial issue: Robert
Liebler (1970), Bruce Cooperstein (1975) and Sylvia Hobart (1987).

A few students
testified that Don gently nudged and engaged them more and more
until they were doing original mathematics. 
One student, Manley Perkel, wrote about his realization during one of Don's advanced courses on permutation groups, that Don was in the midst of working out parts of his new theory on coherent configurations even while simultaneously lecturing on them.Ê Perkel realized only later what a special opportunity Don gave him, as a student, to read early drafts of many of his papers on this subject, especially one component ("Part IV") of the coherent configuration theory.Ê 

Perkel also wrote that ``the mathematics department at the University of Michigan in the '70's was an amazing place to meet distinguished visitors who were there primarily due to Don's influence.Ê Among the long-termÊ visitors during this period were Peter Neumann, Bill Kantor and Peter Cameron, and short-term visitors included Ernie Schult, John Conway, Charles Sims, Bernd Fischer, and Jaap Seidel, among others (including Cheryl Praeger's one-day visit), all there no doubt because of Don.Ê For a graduate student to have the opportunity to meet such distinguished mathematicians was inspiring, to say the least.''   A complete list of Don's visitors would be much longer.

Bob Liebler relates this impression from his grad student years in the late 1960s:
``For many years Higman and Roger  [Lyndon] shared an office on the 3rd floor of Angel hall.  [Lyndon] was a chain smoker and Higman something of an [athlete].  Therefore Higman minimized the
time he spent in the office.  Somehow I discovered that he would hang out in the basement of the
student union often occupying a large table with stacks of papers and computer printouts.  I developed a habit of approaching him there and he was always happy to talk about math or even
politics.''

Don spent a lot of time working
with other mathematicians.
Betty Higman mentions Charles Sims and adds `` \dots
Michio Suzuki spent a lot of time with Don when they were both
working on their Ph.\thinspace  D. at U. of Illinois,
and Hans Zassenhaus was
whom he worked with when he had the National Science Fellowhip
at McGill. He also had contact with Jim Lambek (McGill) and
Jean Miranda [in Montr\'eal].

``Since that time Karl Gruenberg, Bernd Fischer, Christoph Hering,
\dots .  Jaap Seidel was
his main contact at Eindhoven. Most recently he mentioned
Willem Haemers, also from Einhoven, Don Taylor from
Australia was also at Eindhoven when we were there. As I wrote
before I don't have any real idea of how much they influenced
him but I do know he spent a lot of time with each one.

``As far as Jack McLaughlin goes I remember Don saying that they
thought alike and could almost anticipate what the other would
come up with, they were on the same wave length.''

Finally, we mention that Graham Higman and Don Higman both worked on finite groups and combinatorics.  Both had quite heavy eyebrows, a point often remarked on by observers.  Both men had ancestors from the same area in Cornwall, England, though no familial relationship between them has been established.      See \cite{collins}.

\subsubsection{List of Doctoral Students (15 total)}

\ \ \ \ \ James Brooks University of Michigan 1964

Yaotsu Chang University of Michigan 1994 ytchang@isu.edu.tw

Bruce Cooperstein University of Michigan 1975 coop@ucsc.edu

Raymond Czerwinski University of Michigan 1966

David Foulser University of Michigan 1963

Robert Gill University of Michigan 1998

Marshall Hestenes University of Michigan 1967 hestenes@msu.edu

Sylvia Hobart University of Michigan 1987 SAHobart@uwyo.edu

Alan Hoffer University of Michigan 1969

Robert Liebler University of Michigan 1970
liebler@math.colostate.edu

Roger Needham University of Michigan 1992
needham@rio.sci.ccny.cuny.edu

Manley Perkel University of Michigan 1977 manley.perkel@wright.edu

Alyssa Sankey University of Michigan 1992 asankey@unb.ca

Betty Stark University of Michigan 1971 salzberg@ccs.neu.edu

Andrew Vince University of Michigan 1981  vince@math.ufl.edu

\section{Honors}

Donald Higman's academic honors include giving an invited lecture to the 1970
International Congress of Mathematicians in Nice, where he presented his
theory of rank three groups.
Don received the  1975  Alexander von
Humboldt Stiftung Prize. He spent sabbatical and academic leaves in Eindhoven
and Giessen, was a visiting professor at Frankfurt, a visiting senior
scientist at Birmingham and Oxford and a visiting fellow at the Institute for
Advanced Study in Canberra, Australia.

\section{Life in Ann Arbor}

The atmosphere in the University of Michigan Mathematics Department
in the 1950s and 1960s was especially cordial.
There were many parties and frequent expressions of hospitality.

Hirosi Nagao, who visited the University of Michigan in 1958-59, spoke of many
interactions involving Raoul Bott, Donald Higman, Donald Livingstone and Roger
Lyndon. He says ``At  Ann Arbor, Higman, McLaughlin and Raoul Bott used to have
lunch together at the University Union\footnote{`University Union' is currently called the Michigan Union.}. I often joined them and enjoyed
listening their conversations. I was invited many times to [parties] held at
the home of Higman.''

Donald Higman, wife Betty and their five children 
were generous hosts  on many social occasions. Betty worked at the
Library for the Blind and Physically Handicapped. Don had been
an active member of the Flounders and the Ann Arbor Track Club and frequently rode his bicycle.  Don played water polo well into his later years.   Bob Griess remembers many times being at home, noticing Don jogging past as he expelled his breath in loud bursts.  Regular jogging continued for years past Don's 1998 retirement.

\section{Donald Higman publications and preprints}

We found 48 items listed on MathSciNet and 10 additional items.  A historically significant  table [DGH57] is included.  A citation made in this article goes to this section if Donald Higman was an author, and otherwise goes to the Reference section.

\smallskip

 [DGH1] Higman, Donald G. Lattice homomorphisms induced by group homomorphisms. Proc. Amer. Math. Soc. 2, (1951). 467--478. MR0041138 (12,800h) (Reviewer: P. M. Whitman) 20.0X

\smallskip

[DGH2]  Higman, Donald Gordon Focal series in finite groups. Abstract of a thesis, University of Illinois, Urbana, Ill., 1952. ii+1+i pp. 20.0X MR0047644 (13,907c)

\smallskip

[DGH3]  Higman, D. G. Focal series in finite groups. Canadian J. Math. 5, (1953). 477--497. (Reviewer: Graham Higman) 20.0X MR0058597 (15,396b)

\smallskip

[DGH4]  Higman, D. G. Remarks on splitting extensions. Pacific J. Math. 4, (1954). 545--555. (Reviewer: Graham Higman) 20.0X MR0066379 (16,565d)

\smallskip

[DGH5]  Higman, D. G. Modules with a group of operators. Duke Math. J. 21, (1954). 369--376. MR0067895 (16,794b) (Reviewer: B. Eckmann) 20.0X

\smallskip

[DGH6]  Higman, D. G. Indecomposable representations at characteristic $p$. Duke Math. J. 21, (1954). 377--381. MR0067896 (16,794c) (Reviewer: B. Eckmann) 20.0X

\smallskip

[DGH7] Higman, D. G. Induced and produced modules. Canad. J. Math. 7 (1955), 490--508.  MR0087671 (19,390b) (Reviewer: B. Eckmann) 18.0X

\smallskip

[DGH8]Higman, D. G. On orders in separable algebras. Canad. J. Math. 7 (1955), 509--515.  MR0088486 (19,527a) (Reviewer: B. Eckmann) 09.3X

\smallskip

[DGH9]  Higman, D. G. Relative cohomology. Canad. J. Math. 9 (1957), 19--34. MR0083486 (18,715d) (Reviewer: G. P. Hochschild) 09.3X

\smallskip

[DGH10] Higman, D. G.; McLaughlin, J. E. Finiteness of class numbers of representations of algebras over function fields. Michigan Math. J. 6 1959 401--404. MR0109151 (22 \# 39) (Reviewer: W. E. Jenner) 12.00 (16.00)

\smallskip

[DGH11]  Higman, D. G. On isomorphisms of orders. Michigan Math. J. 6 1959 255--257. MR0109174 (22 \# 62)(Reviewer: G. Azumaya) 18.00 (16.00)

\smallskip

[DGH12] Higman, D. G. On representations of orders over Dedekind domains. Canad. J. Math. 12 1960 107--125.  MR0109175 (22 \# 63)(Reviewer: I. Reiner) 18.00 (16.00)

\smallskip

[DGH13]  Higman, D. G.; McLaughlin, J. E. Geometric $ABA$-groups. Illinois J. Math. 5 1961 382--397. MR0131216 (24 \# A1069) (Reviewer: G. Zappa) 50.35

\smallskip

[DGH14]  Higman, D. G. Flag-transitive collineation groups of finite projective spaces. Illinois J. Math. 6 1962 434--446. MR0143098 (26 \# 663) (Reviewer: T. G. Room) 50.60

\smallskip

[DGH15] Higman, Donald G. Finite permutation groups of rank $3$. Math. Z. 86 1964 145--156. MR0186724 (32 \# 4182) (Reviewer: C. Hering) 20.20

\smallskip

[DGH16]  Higman, D. G.; McLaughlin, J. E. Rank $3$ subgroups of finite symplectic and unitary groups. J. Reine Angew. Math. 218 1965 174--189.  MR0175980 (31 \# 256) (Reviewer: P. Dembowski) 20.70

\smallskip

[DGH17]  Higman, Donald G. Primitive rank $3$ groups with a prime subdegree. Math. Z. 91 1966 70--86. MR0218440 (36 \# 1526) (Reviewer: D. L. Barnett) 20.20

\smallskip

[DGH18] Higman, D. G. Intersection matrices for finite permutation groups. J. Algebra 6 1967 22--42. MR0209346 (35 \# 244)  (Reviewer: W. E. Jenner) 20.20

\smallskip

[DGH19]  Higman, Donald G. On finite affine planes of rank $3$. Math. Z. 104 1968 147--149. MR0223971 (36 \# 7018)(Reviewer: H. LŸneburg) 50.60

\smallskip

[DGH20] Higman, Donald G.; Sims, Charles C. A simple group of order $44,352,000$. Math. Z. 105 1968 110--113.  MR0227269 (37 \# 2854) (Reviewer: D. A. Robinson) 20.29

\smallskip

[DGH21]Higman, D. G. Characterization of families of rank 3 permutation groups by the subdegrees. I. Arch. Math. (Basel) 21 1970 151--156.  MR0268260 (42 \#3159) (Reviewer: M. D. Hestenes) 20.20

\smallskip

[DGH22]  Higman, D. G. Characterization of families of rank $3$ permutation groups by the subdegrees. II. Arch. Math. (Basel) 21 1970 353--361. MR0274565 (43 \#328)(Reviewer: M. D. Hestenes) 20.20

\smallskip

[DGH23] Higman, D. G. Coherent configurations. I. Rend. Sem. Mat. Univ. Padova 44 (1970), 1--25.  MR0325420 (48 \#3767)  (Reviewer: Robert A. Liebler) 05B30

\smallskip

[DGH24]  Higman, D. G. Solvability of a class of rank $3$ permutation groups. Nagoya Math. J. 41 1971 89--96. MR0276316 (43 \#2063) (Reviewer: J. P. J. McDermott) 20.20

\smallskip

[DGH25]  Higman, D. G. Characterization of rank $3$ permutation groups by the subdegrees. Representation theory of finite groups and related topics (Proc. Sympos. Pure Math., Vol. XXI, Univ. Wisconsin, Madison, Wis., 1970), pp. 71--72. Amer. Math. Soc., Providence, R.I., 1971. MR0320123 (47 \#8664) (Reviewer: M. D. Hestenes) 20B10

\smallskip

[DGH26] Hestenes, M. D.; Higman, D. G. Rank $3$ groups and strongly regular graphs. Computers in algebra and number theory (Proc. SIAM-AMS Sympos. Appl. Math., New York, 1970), pp. 141--159. SIAM-AMS Proc., Vol. IV, Amer. Math. Soc., Providence, R.I., 1971. MR0340088 (49 \#4844) (Reviewer: Brian Alspach) 05C25

\smallskip

[DGH27]  Higman, D. G. Partial geometries, generalized quadrangles and strongly regular graphs. Atti del Convegno di Geometria Combinatoria e sue Applicazioni (Univ. Perugia, Perugia, 1970), pp. 263--293. Ist. Mat., Univ. Perugia, Perugia, 1971. MR0366698 (51 \#2945) (Reviewer: Jane W. Di Paola) 05B25

\smallskip

[DGH28] Higman, D. G. A survey of some questions and results about rank 3 permutation groups. Actes du Congrs International des MathŽmaticiens (Nice, 1970), Tome 1, pp. 361--365. Gauthier-Villars, Paris, 1971. MR0427435 (55 \#467)  (Reviewer: Shiro Iwasaki) 20B05

\smallskip

[DGH29]  Higman, D. G. Combinatorial considerations about permutation groups. Lectures given in 1971. Mathematical Institute, Oxford University, Oxford, 1972. ii+56 pp. MR0345848 (49 \#10578) (Reviewer: W. M. Kantor) 05B30 (20B25)

\smallskip

[DGH30] Higman, D. G. Remark on Shult's graph extension theorem. Finite groups '72 (Proc. Gainesville Conf., Univ. Florida, Gainesville, Fla., 1972), pp. 80--83. North-Holland Math. Studies, Vol. 7, North-Holland, Amsterdam, 1973.  MR0360346 (50 \#12796) (Reviewer: Chong-Keang Lim) 05C25

\smallskip

[DGH31]  Higman, D. G. Coherent configurations and generalized polygons. Combinatorial mathematics (Proc. Second Australian Conf., Univ. Melbourne, Melbourne, 1973), pp. 1--5. Lecture Notes in Math., Vol. 403, Springer, Berlin, 1974. MR0351864 (50 \#4352) (Reviewer: Brian Alspach) 05B30

\smallskip

[DGH32]  Higman, D. G. Invariant relations, coherent configurations and generalized polygons. Combinatorics (Proc. Advanced Study Inst., Breukelen, 1974), Part 3: Combinatorial group theory, pp. 27--43. Math. Centre Tracts, No. 57, Math. Centrum, Amsterdam, 1974. MR0379244 (52 \#150) (Reviewer: Robert A. Liebler) 05B30

\smallskip

[DGH33]  Higman, D. G. Schur relations for weighted adjacency algebras. Symposia Mathematica, Vol. XIII (Convegno di Gruppi e loro Rappresentazioni, INDAM, Rome, 1972), pp. 467--477. Academic Press, London, 1974. MR0382418 (52 \#3302) (Reviewer: L. Dornhoff) 20C15

\smallskip

[DGH34] Combinatorics. Proceedings of the NATO Advanced Study Institute held at Nijenrode Castle, Breukelen, 8--20 July 1974. Edited by M. Hall, Jr. and J. H. van Lint. NATO Advanced Study Institutes Series, Series C: Mathematical and Physical Sciences, Vol. 16. D. Reidel Publishing Co., Dordrecht-Boston, Mass.; Mathematical Centre, Amsterdam, 1974. viii+482 pp. 05-06 (94A20)  MR0387065 (52 \#7912)

\smallskip

[DGH35]  Higman, D. G. Coherent configurations. I. Ordinary representation theory. Geometriae Dedicata 4 (1975), no. 1, 1--32. MR0398868 (53 \#2719) (Reviewer: Robert A. Liebler) 05B99

\smallskip

[DGH36]  Higman, D. G. A monomial character of Fischer's baby monster. Proceedings of the Conference on Finite Groups (Univ. Utah, Park City, Utah, 1975), pp. 277--283. Academic Press, New York, 1976.  MR0409627 (53 \#13379) (Reviewer: Bhama Srinivasan) 20D05

\smallskip

[DGH37]  Higman, D. G. Coherent configurations. II. Weights. Geometriae Dedicata 5 (1976), no. 4, 413--424.  MR0437368 (55 \#10299) (Reviewer: Robert A. Liebler) 05B99 (20C99)

\smallskip

[DGH38]  Higman, D. G. Lectures on permutation representations. Notes taken by Wolfgang Hauptmann. Vorlesungen aus dem Mathematischen Institut Giessen, Heft 4. Mathematisches Institut Giessen, Giessen, 1977. 106 pp. MR0470042 (57 \#9810) (Reviewer: Michael Klemm) 20BXX

\smallskip

[DGH39] Higman, D. G. Systems of configurations. Proceedings, Bicentennial Congress Wiskundig Genootschap (Vrije Univ., Amsterdam, 1978), Part I, pp. 205--212, Math. Centre Tracts, 100, Math. Centrum, Amsterdam, 1979. 20C20 (05B99 05C20)  MR0541394 (80m:20008)

	Note: [DGH39] Systems of Configurations is the title of two preprints, which differ in their covers, but appear to be the same.

\smallskip

[DGH40]  Higman, D. G. Admissible graphs. Finite geometries (Pullman, Wash., 1981), pp. 211--222, Lecture Notes in Pure and Appl. Math., 82, Dekker, New York, 1983.  MR0690807 (84c:05026) (Reviewer: Francis Buekenhout) 05B25 (05C25 51A10)

\smallskip

[DGH41] Higman, D. G. Coherent algebras. Linear Algebra Appl. 93 (1987), 209--239. 15A30 	MR0898557 (89d:15001)

\smallskip

[DGH42]  Higman, D. G. Strongly regular designs and coherent configurations of type $[{3\atop {\;}}\;\;{2\atop 3}]$. European J. Combin. 9 (1988), no. 4, 411--422. MR0950061 (89i:05070) (Reviewer: Dieter Jungnickel) 05B30

\smallskip

[DGH43]  Haemers, W. H.; Higman, D. G. Strongly regular graphs with strongly regular decomposition. Linear Algebra Appl. 114/115 (1989), 379--398. MR0986885 (90b:05108) (Reviewer: Dieter Jungnickel) 05C75 (05B30)

\smallskip

[DGH44]  Higman, D. G. Computations related to coherent configurations. Proceedings of the Nineteenth Manitoba Conference on Numerical Mathematics and Computing (Winnipeg, MB, 1989). Congr. Numer. 75 (1990), 9--20. MR1069158 (92a:05032) (Reviewer: Arnold Neumaier) 05B30 (05E30 15A30 20C99)

\smallskip

[DGH45] Higman, D. G. Weights and $t$-graphs. Algebra, groups and geometry. Bull. Soc. Math. Belg. SŽr. A 42 (1990), no. 3, 501--521. MR1316208 (96g:20006)  (Reviewer: Hugo S. Sun) 20C15 (05C25)

\smallskip

[DGH46]  Haemers, W. H.; Higman, D. G.; Hobart, S. A. Strongly regular graphs induced by polarities of symmetric designs. Advances in finite geometries and designs (Chelwood Gate, 1990), 163--168, Oxford Sci. Publ., Oxford Univ. Press, New York, 1991.  MR1138741 (92i:05028) (Reviewer: Dina Smit-Ghinelli) 05B05 (05E30)

\smallskip

[DGH47] Higman, D. G. Rank $5$ association schemes and triality. Linear Algebra Appl. 226/228 (1995), 197--222.   MR1344562 (96j:05117) (Reviewer: Sanpei Kageyama) 05E30

\smallskip

[DGH48]  Higman, D. G. Strongly regular designs of the second kind. European J. Combin. 16 (1995), no. 5, 479--490.  MR1345694 (96g:05035) (Reviewer: Esther R. Lamken) 05B30

\smallskip

\smallskip

{\bf  The items below are not referred to in MathScinet.  Except for [DGH50], they are not published as far as we know.   The first three have dates and the others may not.  Some physical characteristics are noted.  }

\smallskip

[DGH49] PREPRINT, Part {$\underline { \ \ \ }$}, Homogeneous Configurations of Rank 4, [Don Higman told Manley Perkel that this was supposed to be part 4 of a series] 36 pages (hand written) . Don Higman gave this to Manley Perkel in 1974 to proofread. It is related to notes from a course in 1972-73 at University of Michigan.

\smallskip

[DGH50] REPRINT, Monomial Representations; this was published
in Finite Groups, Sapporo and Kyoto, 1974. Proceedings of the Taniguchi
International Symposium, edited by N, Iwahori. copyright 1976, Japan Society
for the Promotion of Science (Japanese NSF)

\smallskip

[DGH51] PREPRINT, DATED 8/26/97 Some highly symmetric chamber systems, 8 pages.
	from Don Higman's collection, obtained by Bob Griess

\smallskip

[DGH52] PREPRINT, A note on regular 3-graphs, 7 pages (typed) 	contributed by Alyssa Sankey

\smallskip

[DGH53] PREPRINT, The parabolics of a semi-coherent configuration, 18 pages (seems to be dot matrix printed) 	 contributed by Alyssa Sankey

\smallskip

[DGH54] PREPRINT, Relation configurations and relation algebras, 27 pages (typed -- signature under title) 	contributed by Alyssa Sankey

\smallskip

[DGH55] PREPRINT, Uniform association schemes, 23 pages (word processed?) 	contributed by Alyssa Sankey

\smallskip

[DGH56] PREPRINT, Untitled, 16 pages (the table of contents starts with color schemes and morphisms, ends with homology and weights)(word processed, with handwriting) 	contributed by Alyssa Sankey

\smallskip

[DGH57]
Title: THIS IS A LIST OF THE KNOWN RANK THREE GROUPS FOR DEGREES UP TO 10000.      A copy was provided by Francis  Buekenhout, who gives the date 17 July, 1968.  The print had faded and was difficult to read.  A transcription of it, done in 2007, consisting of  most of the information on the table, is on the web site of Robert Griess.  We thank Ching Hung Lam for arranging the typing.

\smallskip

[DGH58] A table of parameter values for rank 3 groups;
In the mid 1990s,
Don Higman gave a table of parameter values for rank 3 groups to Bob Griess for inclusion in his book \cite{Griess12}, page 125.  This is a smaller set of data than [DGH57].


\end{document}